\documentclass[preprint,3p,11pt,sort&compress]{elsarticle}

 \usepackage{hyperref}
 \hypersetup{backref, colorlinks=true, linkcolor=blue}

\usepackage{stfloats}
\usepackage{subcaption}

 \usepackage{pict2e}

\usepackage{graphicx} 
\usepackage{booktabs} 
\usepackage{microtype} 
\usepackage{amssymb,amsmath}
\usepackage[ruled,vlined]{algorithm2e}

\begin{document}
\begin{center}
{\bf \Large Hybrid Weight Window Techniques for Time-Dependent \\
\vspace{0.1cm}
 Monte Carlo Neutronics}
\end{center}

\author[ncsu]{Caleb S. Shaw}
\ead{cashaw4@ncsu.edu}
\author[ncsu]{Dmitriy Y. Anistratov}
\ead{anistratov@ncsu.edu}

\address[ncsu]{Department of Nuclear Engineering,
North Carolina State University Raleigh, NC}

\begin{frontmatter}
\begin{abstract}
Efficient variance reduction of Monte Carlo simulations is desirable to avoid wasting computational resources. This paper presents an automated weight window algorithm for solving time-dependent particle transport problems. The weight window centers are defined by a hybrid forward solution of the discretized low-order second moment (LOSM) problem. The second-moment (SM) functionals defining the closure for the LOSM  equations are computed by Monte Carlo solution. A filtering algorithm is applied to reduce noise in the SM functionals. The LOSM equations are discretized with first- and second-order time integration methods. We present numerical results of the AZURV1 benchmark. The hybrid weight windows lead to a uniform distribution of Monte Carlo particles in space. This causes a more accurate resolution of wave fronts and regions with relatively low flux.
\end{abstract}
\begin{keyword}
dynamic Monte carlo, 
particle transport, 
variance reduction, 
weight windows, 
second moment method, 
hybrid techniques
\end{keyword}
\end{frontmatter}

\section{Introduction}\label{sec:intro}

 Monte Carlo methods are used for solving time-dependent neutron transport problems due to their many advantages, such as large time steps, continuous energy physics and lack of discretization error. However, these methods converge relatively slowly with increasing number of particles. Analog Monte Carlo distributes particles inefficiently such that regions with sources end up oversampled, leaving other regions undersampled. This drives the development of variance reduction techniques such as implicit capture and weight windows  \cite{wagner-haghighat,wagner-peplow-mosher-2014}. Both of these techniques alter the weight of the Monte Carlo particles. The distribution of the Monte Carlo particles in phase-space can then be modified to improve efficiency. In this paper, we focus on the definition of 
weight windows using an auxiliary forward solution. This technique can obtain a uniform particle distribution and variance. It has also been shown that weight window centers can be efficiently defined by the numerical solution of hybrid moment equations such as the quasi-diffusion (aka variable Eddington factor)  low-order equations
 \cite{Cooper_Larsen_2001,Cooper_diss_1999,Wollaber_2008,gol'din-cmmp-1964,auer-mihalas-1970}.

We propose applying a different low-order problem that is used to formulate the second moment (SM) method \cite{smm-1976,Cefus_Larsen_1989}. The differential operator of the low-order SM (LOSM) equations is self-adjoint. This feature has advantages for solving multi-dimensional particle transport problems. To approximate LOSM equations in time we apply two time-integration scheme, Backward Euler, and  Crank-Nicolson. A second-order finite volume scheme is used for spatial discretization of the LOSM equations. To reduce effects of noise on the 
hybrid LOSM problem, we apply a noise filtering technique
to coefficients of SM closure terms and initial conditions at a timestep. The proposed automated weight windows 
defined by the hybrid solution of the LOSM equations improves the cell wise figure of merit. The developed automatic weight window algorithms have been implemented in the open source python-based dynamic neutron transport code MC/DC \cite{Joss}.

The remainder of the paper is organized as follows.
In Section \ref{sec:methodology} we introduce the methodology used with a focus on the second moment method and the weight window definition. Section \ref{sec:results} presents the numerical results for the chosen test problem where we discuss the performance of the algorithm. Section \ref{sec:conclusion} concludes with a summary of our findings and potential directions for future work.
 
\section{Methodology} \label{sec:methodology}

\subsection{Second Moment Method}
 
The time dependent behavior of the ensemble of particles is described by the neutron transport equation. Here it is in 1D 1-group form with isotropic scattering and isotropic source:
 \begin{equation}  
    \frac{1}{v}\frac{\partial\psi}{\partial t}(x,\mu,t)+\mu\frac{\partial \psi}{\partial x}(x,\mu,t)+\Sigma_t(x,t)\psi(x,\mu,t)=\frac{1}{2}\big((\Sigma_s(x,t)+ 
    \nu_f \Sigma_f(x,t)  )\phi(x,t) +q(x,t)\big),
    \label{eqn:HO_transport}
\end{equation}
\[
x\in[0,X],\text{ } \mu\in[-1,1],\text{ } t\geq 0,
\]
and boundary conditions (BCs)
\begin{equation*}
    \psi(0,\mu,t) = \psi_L(\mu,t),\text{ for } \mu> 0, \quad
    \psi(X,\mu,t) = \psi_R(\mu,t),\text{ for } \mu< 0, 
\end{equation*}
and initial conditions (ICs)
\begin{equation*}
    \psi(x,\mu,0)  = \psi_0(x,\mu),
\end{equation*}
where $\psi$ is the particle angular flux, $x$ is the particle spatial position,
$\mu$ is the directional cosine of particle motion, $t$ is time, $v$ is the particle speed,
$\Sigma_t$, $\Sigma_s$, $\Sigma_f$ are total, scattering, and fission cross sections, respectively,
$\nu_f$ is the expected number of neutrons born in fission,
$q$ is the external source of particles. The low-order equations are derived by integrating the transport equation (Eq. \ref{eqn:HO_transport})
with weights 1 and $\mu$ over $-1 \le \mu \le 1$.
The  LOSM  equations for the scalar flux, $\phi(x,t) = \int_{-1}^{1} \psi (x,\mu,t) d \mu$, and current, $J(x,t) = \int_{-1}^{1} \mu \psi (x,\mu,t) d \mu$, 
 are given by \cite{smm-1976}
\begin{equation}
    \frac{1}{v}\frac{\partial \phi}{\partial t}(x,t) + \frac{\partial J}{\partial x}(x,t) + \left(\Sigma_t(x,t)-\Sigma_s(x,t)-\nu_f\Sigma_f(x,t)\right)\phi(x,t) = q(x,t),
    \label{eqn:balance}
\end{equation}
\begin{equation}
    \frac{1}{v}\frac{\partial J}{\partial t}(x,t) + \frac{1}{3}\frac{\partial \phi}{\partial x}(x,t) +\Sigma_t(x,t)J(x,t) = \frac{\partial F}{\partial x}(x,t),
    \label{eqn:moment}
\end{equation}
with the associated BCs
\begin{equation}
    J(0,t) = -\frac{1}{2}\phi(0,t)+2J_L(t)+P_L(t),\quad J(X,t) = \frac{1}{2}\phi(0,t)+2J_R(t)-P_R(t),
     \label{losm-bcs}
\end{equation}
and ICs
\begin{equation*} 
    \phi(x,0) = \phi_0(x), \quad J(x,0) = J_0(x). 
    \label{losm-ics}
\end{equation*}
The terms with the second angular moment of the high-order solution are closed with the function defined by
\begin{equation} \label{F}
   F(x,t) = \int_{-1}^{1}\left(\frac{1}{3}-\mu^2\right)\psi(x,\mu,t)d\mu \, .
\end{equation}
The BCs (Eq. \ref{losm-bcs}) are formulated with 
\begin{equation*}
   P_L(t) = \int_{-1}^{1}\left(\frac{1}{2}-|\mu|\right)\psi(0,\mu,t)d\mu,\quad
   P_R(t) = \int_{-1}^{1}\left(\frac{1}{2}-|\mu|\right)\psi(X,\mu,t)d\mu \, .
\end{equation*}

The LOSM equations are approximated in time using two different methods. We apply the  Backward Euler (BE) time-integration scheme. Its stability properties are not dependent on the timestep size. The LOSM equations discretized with the BE scheme are given by
\begin{equation}
    \frac{1}{v\Delta t^n}\left(\phi^n-\phi^{n-1}\right) + \frac{dJ^n}{dx}+\left(\Sigma_t^n- \Sigma_s^n-\nu_f^n\Sigma_f^n \right)\phi^n=q^n,
  \label{eqn:LOSM-BE-1}
\end{equation}
\begin{equation}
    \frac{1}{v\Delta t^n}\left(J^n-J^{n-1}\right) +\frac{1}{3}\frac{d\phi^n}{dx}+\Sigma_t^nJ^n=\frac{dF^n}{dx},
  \label{eqn:LOSM-BE-2}
\end{equation}
where $n$ is the index of the time layer and e.g. $\phi^n(x) =  \phi^n(x,t^n)$. Another temporal scheme for the LOSM equation is derived using a second-order Crank-Nicolson (CN) method. This scheme reduces effects of the numerical diffusion compared to the BE method and improves the solution accuracy for the same value of the timestep. The LOSM equations approximated by means of the CN scheme have the following form:
\begin{multline}
    \frac{1}{v\Delta t^n}\left(\phi^n-\phi^{n-1}\right)+\frac{1}{2}\left[\frac{dJ^n}{dx}+(\Sigma_t^n-\Sigma_s^n-\nu_f^n\Sigma_f^n)\phi^n\right] =\\
    \frac{1}{2}(q^n+q^{n-1})-\frac{1}{2}\left[\frac{dJ^{n-1}}{dx}+(\Sigma_t^{n-1}-\Sigma_s^{n-1}-\nu_f^{n-1}\Sigma_f^{n-1} )\phi^{n-1}\right],
\label{eqn:LOSM-CN-1}
\end{multline}
\begin{equation}
    \frac{1}{v\Delta t^n}\left(J^n-J^{n-1}\right)+\frac{1}{2}\left[\frac{1}{3}\frac{d\phi^n}{dx}+\Sigma_t^nJ^n\right] = \frac{1}{2}\left[\frac{dF^n}{dx}+\frac{dF^{n-1}}{dx}\right]-\frac{1}{2}\left[\frac{1}{3}\frac{d\phi^{n-1}}{dx}+\Sigma_t^{n-1}J^{n-1}\right].
\label{eqn:LOSM-CN-2}
\end{equation}

To spatially discretize the LOSM equations we use a second order finite volume scheme. We define a spatial mesh $\{x_{i-1/2}\}_{i=0}^I$. The semi-discretized balance equation (Eqs. \eqref{eqn:LOSM-BE-1}  and \eqref{eqn:LOSM-CN-1}) is integrated over the $i^{th}$ spatial $[x_{i-1/2}, x_{i+1/2}]$. The semi-discretized first-moment equation (Eqs. \eqref{eqn:LOSM-BE-2} and \eqref{eqn:LOSM-CN-2}) is integrated over halves of spatial cells. The resulting discretized equations are defined for cell averaged fluxes, $\phi_i^n~=~\frac{1}{\Delta x_i} \int_{x_{i-1/2}}^{x_{i+1/2}}\phi^n(x)dx$, and cell-edge currents, $J_{i\pm1/2}^n= J^n( x_{i\pm1/2} )$, and are given by 
\begin{equation} 
    J_{i+1/2}^n-J_{i-1/2}^n+(\hat{\Sigma}_{t,i}^n-\Sigma_{s,i}^n- \nu_{f,i}^n\Sigma_{f,i}^n )\Delta x_i\phi_i^n=Q_{0,i}^n\Delta x_i \, ,
 \label{FV-LOSM-1}
\end{equation}
\begin{equation}
\frac{1}{3} \Big( \phi_{i}^n - \phi_{i-1}^n\Big)   +   \hat{\Sigma}_{t,i-1/2}^n \Delta x_{i-1/2} J_{i-1/2}^n =
Q_{1,i-1/2}^{n} \, ,
 \label{FV-LOSM-2}
\end{equation} 
\begin{equation*}
\hat{\Sigma}_{t,i}^n= \begin{cases}
    \Sigma_{t,i}^n +\frac{1}{v\Delta t^n},\quad \text{for BE},\\
     \Sigma_{t,i}^n +\frac{2}{v\Delta t^n},\quad \text{for CN},\\
\end{cases}
\quad
  \hat \Sigma_{t,i-1/2}^n =\frac{ \hat \Sigma_{t,i-1}^n \Delta x_{i-1} + \hat \Sigma_{t,i}^n\Delta x_i }{\Delta x_{i-1} + \Delta x_i} \, , 
\end{equation*}
\begin{equation*}
\Delta x_{i} =    x_{i+1/2} - x_{i-1/2} \, , \quad
\Delta x_{i-1/2} = \frac{1}{2}( \Delta x_{i+1} +  \Delta x_i) \, ,  \quad i=1,...,N_x \, ,
\end{equation*}
and the right-hand sides are defined by
\begin{align*}
    Q&_{0,i}^n = \begin{cases}
        q_i^n+\frac{\phi_i^{n-1}}{v\Delta t^n},\quad \text{for BE,}\\
         q_i^n+q_i^{n-1} + 
        \frac{1}{\Delta x_i}(J_{i-1/2}^{n-1}-J^{n-1}_{i+1/2}) + \Big(\Sigma_{s,i}^{n-1}+\nu_{f,i}^{n-1}\Sigma_{f,i}^{n-1}+
        \frac{2}{v\Delta t^n}-\Sigma_{t,i}^{n-1}\Big)\phi_{i}^{n-1}, 
        \quad \text{for CN},
    \end{cases}
\end{align*}
\begin{align*}
    Q&_{1,i-1/2}^n=
    &\begin{cases}
        \frac{J_{i-1/2}^{n-1}}{v\Delta t^n},\quad \text{for BE,}\\
        \frac{1}{\Delta x_{i-1/2}} \! \Big( \! F_i^{n-1} \! - \! F_{i-1}^{n-1} \! + \! \frac{1}{3}(\phi_{i-1}^{n-1} \! - \!  \phi_i^{n-1}) \! \Big)  \! + \! 
        \Big( \! \frac{2}{v\Delta t^n} \! - \! \frac{\Sigma_{t,i-1}^{n-1}\Delta x_{i-1}+\Sigma^{n-1}_{t,i}\Delta x_i}{\Delta x_{i-1}+\Delta x_i} \! \Big)J_{i-1/2}^{n-1},\ \text{for CN}, 
    \end{cases}
\end{align*}
where
\begin{equation} 
F_i^n= \frac{1}{\Delta x_i} \int_{x_{i-1/2}}^{x_{i+1/2}}F^n(x)dx \, .
\label{F_i}
\end{equation}

The space-time integration scheme applied to discretize the LOSM equations defines the grid functions of the solution and closure terms in spatial cells and at time layers. At each timestep, these key quantities are computed by Monte Carlo algorithms using an ensemble of particle histories. The scalar flux $\phi_i^{n-1}$, and the functionals for closure terms $F_i^{n}$ and $F_i^{n-1}$ are cell averages computed at discrete temporal points rather than time averages. The current $J_{i\pm1/2}^{n-1}$ is similarly computed along the cell edges at the time layer.

Monte Carlo simulations use track-length tallies to numerically integrate quantities with greater statistical reliability than census or crossing tallies. Track length tallies capture information by accumulating data along particle paths within cells, whereas census and crossing tallies rely on particles crossing a temporal boundary or spatial face. However track-length tallies produce time-averaged quantities, which do not directly provide the values at discrete time points for the LOSM equations. There exist numerical techniques which use a second order backward difference (BDF-2) formula to extrapolate quantities obtained from Monte Carlo photon transport tallies \cite{mcclarren2013temperature}.

In the developed hybrid algorithm, we use the particle locations at the time layer to compute cell averaged
 $\phi_i^n$ and $F_i^n$. This is a direct Monte Carlo estimation of these desired quantities. However evaluating $J_{i\pm1/2}^n$ is less straightforward, since particles being on the cell edge at the time layer is exceedingly rare. Instead, we tally $J_i^n= \frac{1}{\Delta x_i} \int_{x_{i-1/2}}^{x_{i+1/2}}J^n(x)dx$, the spatially averaged current on the time layer, and use linear interpolation to calculate $J_{i\pm1/2}^n$. 
  
  It is important to note that $F_i^{n}$ is not immediately available from the Monte Carlo histories of the previous timestep. Consequently, at the beginning of calculations on the $n^{th}$ timestep, we lag this quantity  
and compute weight windows with the solution of the resulting semi-implicit time approximation of LOSM equations. After half the histories have been run during the timestep, $F_i^{n}$ is calculated using the Monte Carlo solution on the current timestep.
Then we compute the hybrid solution by solving fully-implicit discretized LOSM equations and update weight windows. 

The parameters defining the hybrid LOSM problems computed by Monte Carlo techniques 
 inherently have statistical noise. This affects the numerical solution of the LOSM equations.
To reduce high frequency noise in coefficients of the LOSM equations and initial data at a time step,  a moving average filter is applied to  
 $\phi_i^{n-1}$, $J_{i\pm1/2}^{n-1}$, and $F_i^{n}$,
which are  all  calculated  by Monte Carlo algorithms. Any of these grid functions is represented as a discrete vector, $\mathbf{f}=\{f_i\}_{i=1,...,I}$. For each element, $f_i$, the filtered value, $\tilde{f}_i$ is given by \cite{Smith-book}
\begin{equation}
	\label{eqn:moving_average}
	\tilde{f}_i = \frac{1}{2k+1}\sum_{j=-k}^k f_{i-j} \, ,
	\vspace{0.25cm}
\end{equation}
where $k$ is the filter width. On the edges, the window would extend outside of the range of the vector. The width is adaptively shrunk to avoid this issue. We note that  while smoothing may bias the auxiliary hybrid solution, since this solution will only be used to set the weight windows, the final Monte Carlo solution will remain unbiased.

\subsection{Automated Weight Window Methodology}
In this section we will describe how to construct the weight windows based on
the hybrid forward solution  \cite{Cooper_Larsen_2001}. We begin by imposing a spatial mesh, and in each cell we define the following quantities.
\vspace{-0.15cm}
\begin{equation}
    ww_i^n = \bigg[\frac{\tilde{\phi}_i^n}{\max\limits_i \tilde{\phi}_i^n}\bigg]
    \times(1-\varepsilon_{min})+\varepsilon_{min}= \text{ center of the weight window in cell $i$ on timestep $n$,}
\end{equation}
\begin{equation}
    ww_i^{n,ceiling} = ww_i^n\times \rho = \text{ ceiling of the weight window in cell $i$ on timestep $n$,}
\end{equation}
\begin{equation}
    ww_i^{n,floor} = \frac{ww_i^n}{ \rho }= \text{ floor of the weight window in cell $i$ on timestep $n$,}
\end{equation}
where 
\begin{equation*}
   \tilde{\phi}_i^n = \text{ an approximation of the scalar flux in cell $i$ on timestep $n$},
\end{equation*}
\begin{equation*}
    \rho = \text{ the window width parameter} \geq 1,
\end{equation*}
\begin{equation*}
    \varepsilon_{min} = \text{ the minimum center parameter} < 1.
\end{equation*}

Whenever a particle enters a new cell, its weight is compared to $ww_i^{n,floor} $ and $ww_i^{n,ceiling}$. If its weight is larger than $ww_i^{n,ceiling}$, the particle is split such that the daughter particles weights are within the window. If the particle weight is lower than $ww_i^{n,floor}$, the particle will undergo a rouletting process that will either kill it, or increase its weight to be within the window. Both splitting and rouletting preserve the total weight and do not introduce bias into the solution.

 The $\varepsilon_{min}$ parameter is introduced to avoid window centers being zero or near zero which would result in very small weight windows. This leads to excessive particle splitting which can be computationally infeasible. Using the minimum center modification allows the user to control the minimum window centers, while preserving the shape from the auxiliary solution. If $\varepsilon_{min}$ is too small, or zero, it can lead to inefficiencies and in some cases failure due to overwhelming particle splitting.

 In regions of the problem with significant flux gradients, such as wave fronts, the window centers can drop orders of magnitude over a couple cells. This is undesirable since it can again cause excessive splitting leading to computational inefficiency. To prevent very sharp drop offs in $ww_{i}^{n}$, 
 we can use the modification developed for Marshak waves, which artificially raises the window centers at the front to reduce splitting.    
 The modified window center in cell $i$ on timestep  $n$ is given by \cite{Wollaber_2008} 
 \begin{equation}
    \hat{ww}_{i}^{n}=ww_{i}^{n}\times\left[1 + \left(\frac{1}{\varepsilon}-1\right)e^{-(ww_{i}^{n}-w_{min})/\varepsilon}\right],
 \vspace{-0.2cm}
 \end{equation}
where  $\varepsilon$ and $w_{min}$ are  the parameters of the modification.
 
In the rest of this paper we will refer to the Monte Carlo solution computed using weight windows defined by the hybrid solution as LOSM-BE when using Backward Euler and as LOSM-CN when using Crank-Nicolson. We will compare these methods to the established approach of using the Monte Carlo solution from the previous timestep to define the weight windows \cite{Landman_McClarren_Madsen_Long_2014}. For the problem presented below, we also compare with the weight windows defined by the analytic solution. We will refer to the Monte Carlo solution computed using these two methods as WW-Previous and WW-Analytic respectively.

\section{Numerical Results}\label{sec:results}

The AZURV1 suite of benchmarks is often used to test dynamic neutronics codes since it contains time-dependent infinite medium problems with an analytic solutions \cite{azurv1}. We have chosen a version from this suite that is driven by an impulse at $t=0$, meaning the only mechanism of producing particles for $t>0$ is multiplication due to fission. This problem is supercritical and hence population control techniques must be used to avoid a runaway particle population. We take advantage of the analytic solution to define weight windows without discretization error or statistical noise for comparison with the hybrid-LOSM generated windows. We also use the analytic solution to compute the error in the hybrid solution. 

The spatial domain $(-20.5 \, cm\leq x \leq 20.5 \, cm)$ contains one material.
The tally mesh is uniform with   201 cells. The material properties are $\Sigma_s = \frac{1}{3}$ cm$^{-1}$, $\Sigma_f = \frac{1}{3}$ cm$^{-1}$, $\Sigma_t = 1$ cm$^{-1}$, and $\nu_f = 2.3$. This means the effective scattering ratio is $c=1.1$. 
The speed of particles is $v =1$  cm~s$^{-1}$. The problem is solved over $0 \, s\leq t\leq 10 \, s$ and there are 20 timesteps with $\Delta t = 0.5 \, s$. The problem is driven by a point source at $x=0$ and time $t=0$ with strength 1 cm$^{-3}$ s$^{-1}$. The boundary conditions are reflective on both sides. The weight window mesh is the same as the tally mesh. The parameters used for the weight window center modifications were $\varepsilon_{min} = 10^{-4}$, $\varepsilon = 10^{-4}$, and $w_{min}= 10^{-4}$. The moving average filter was applied with $k=2$.

\begin{figure}[t!]
\centering
\begin{minipage}{0.5\textwidth}
    \centering
    \includegraphics[width=\linewidth]{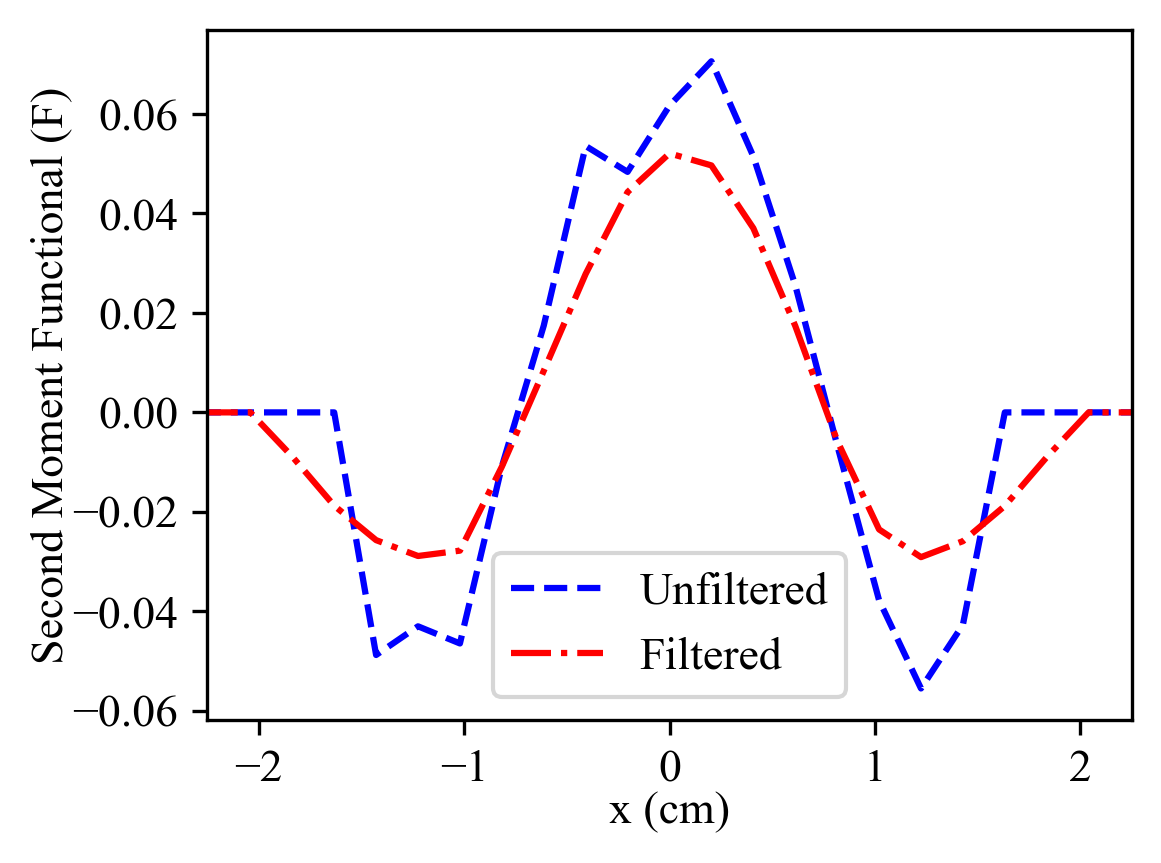}
    \caption{Second moment functional $F^n$ at $t = 2$ s}
    \label{fig:sm_factor3}
\end{minipage}\hfill
\begin{minipage}{0.5\textwidth}
    \centering
    \includegraphics[width=\linewidth]{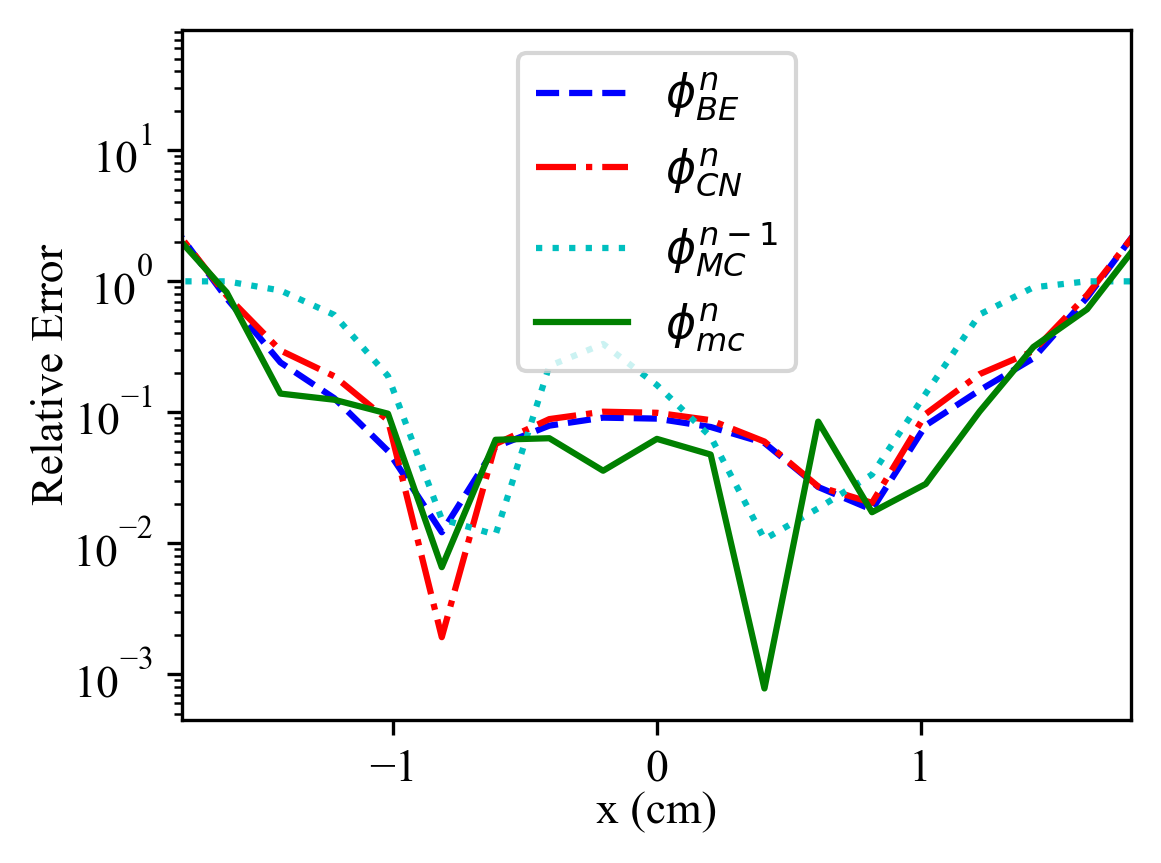}
    \caption{Solution relative error $t = 2$ s}
    \label{fig:det_flux3}
\end{minipage}\hfill
\begin{minipage}{0.5\textwidth}
    \centering
    \includegraphics[width=\linewidth]{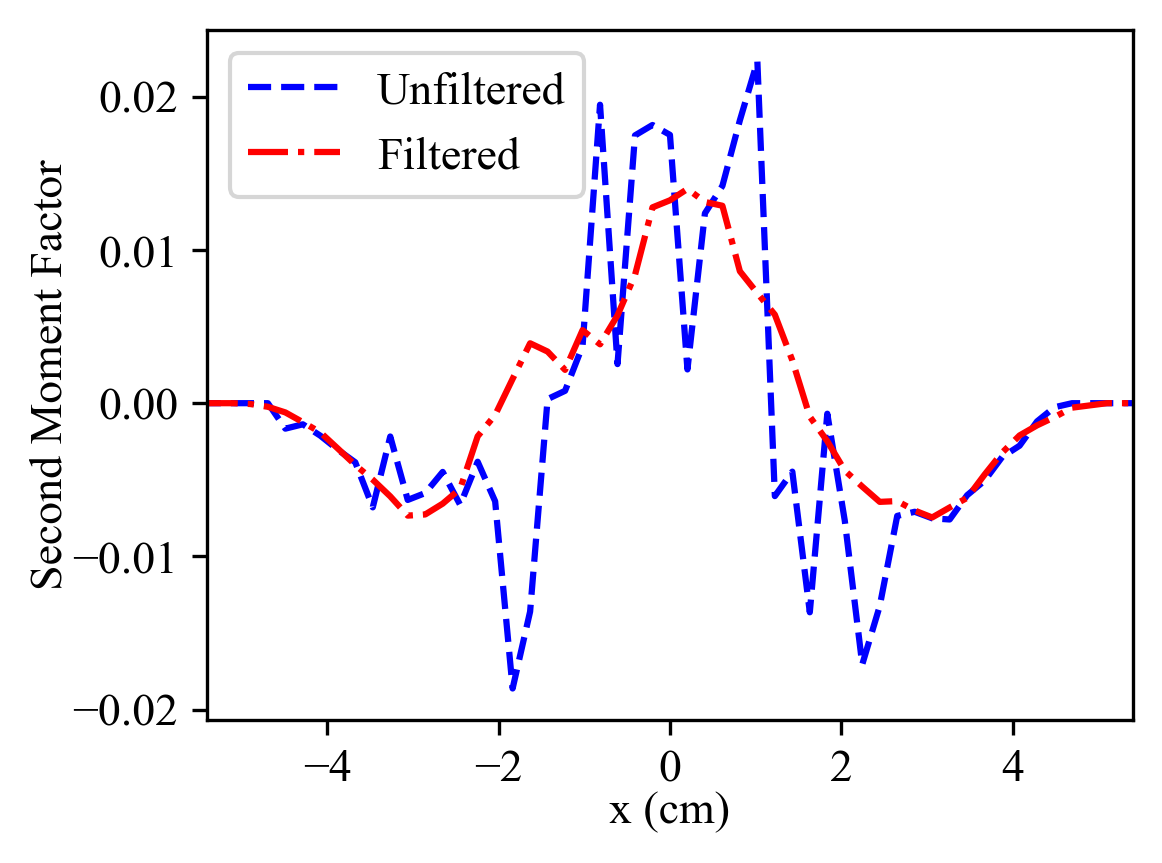}
    \caption{Second moment functional $F^n$ at $t = 5$ s}
    \label{fig:sm_factor9}
\end{minipage}\hfill
\begin{minipage}{0.5\textwidth}
    \centering
    \includegraphics[width=\linewidth]{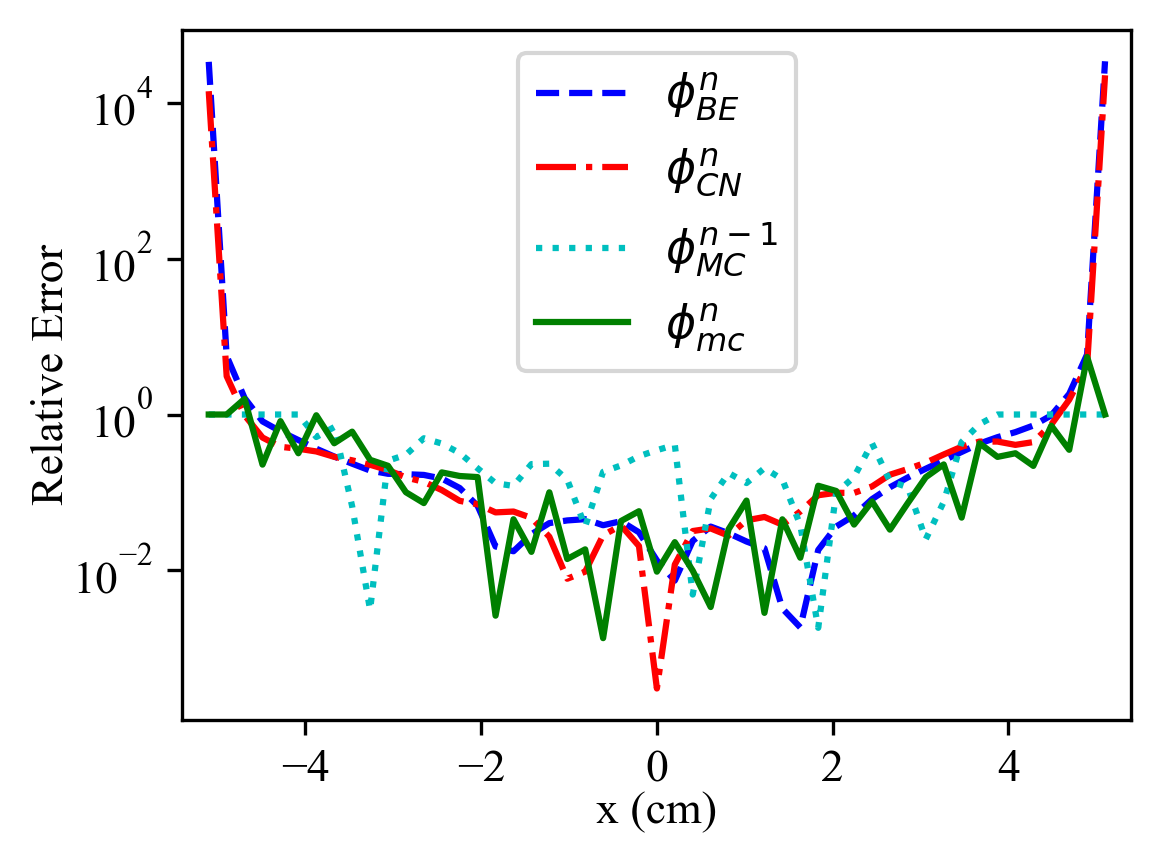}
    \caption{Solution relative error $t = 5$ s}
    \label{fig:det_flux9}
\end{minipage}\hfill
\begin{minipage}{0.5\textwidth}
    \centering
    \includegraphics[width=\linewidth]{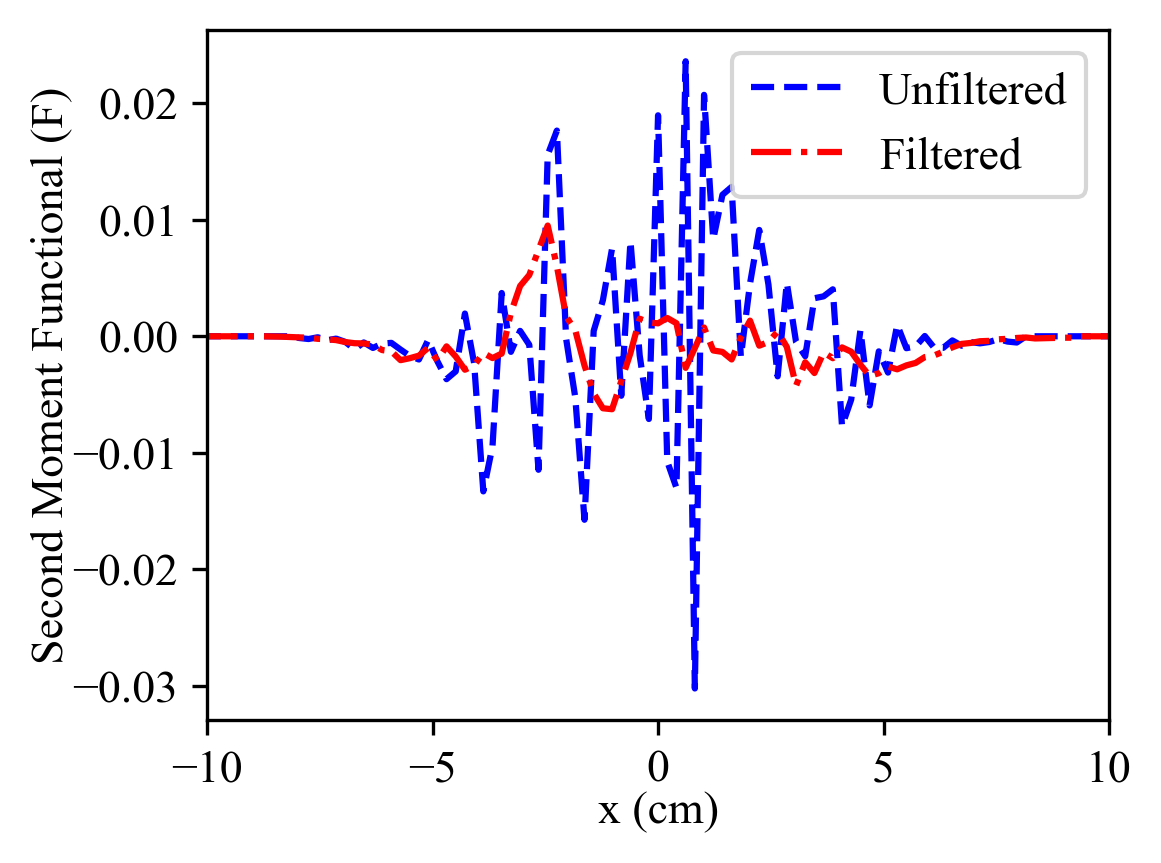}
    \caption{Second moment functional $F^n$ at $t = 10$ s}
    \label{fig:sm_factor-1}
\end{minipage}\hfill
\begin{minipage}{0.5\textwidth}
    \centering
    \includegraphics[width=\linewidth]{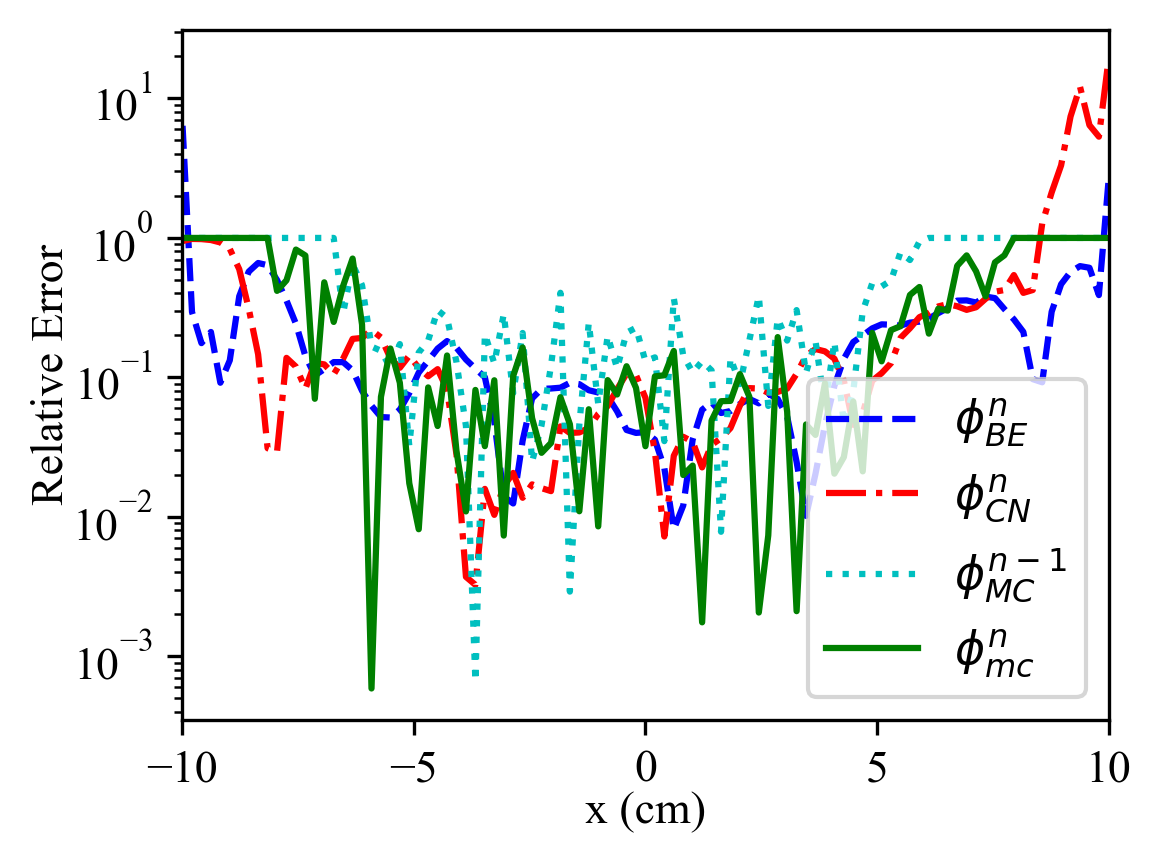}
    \caption{Solution relative error $t = 10$ s}
    \label{fig:det_flux-1}
\end{minipage}\hfill
\end{figure}

In Figures \ref{fig:sm_factor3}, \ref{fig:sm_factor9}, and \ref{fig:sm_factor-1} the effects of the moving average filter on the second moment functional are shown, demonstrating a reduction of noise while preserving the underlying structure. Figures \ref{fig:det_flux3}, \ref{fig:det_flux9}, and \ref{fig:det_flux-1} show the relative error in hybrid LOSM-BE and LOSM-CN solutions   ($\phi_{BE}^{n}$, $\phi_{CN}^{n}$) and Monte Carlo   solutions  ($\phi_{MC}^{n-1}$), used as auxiliary solution $\tilde \phi^n$ to define weight windows. These plots indicate that the hybrid solution is of similar accuracy to the final Monte Carlo solution on the timestep, using information from only half of the histories.

\begin{figure}[t!]
    \centering
    \begin{minipage}{0.5\textwidth}
        \centering
        \includegraphics[width=\linewidth]{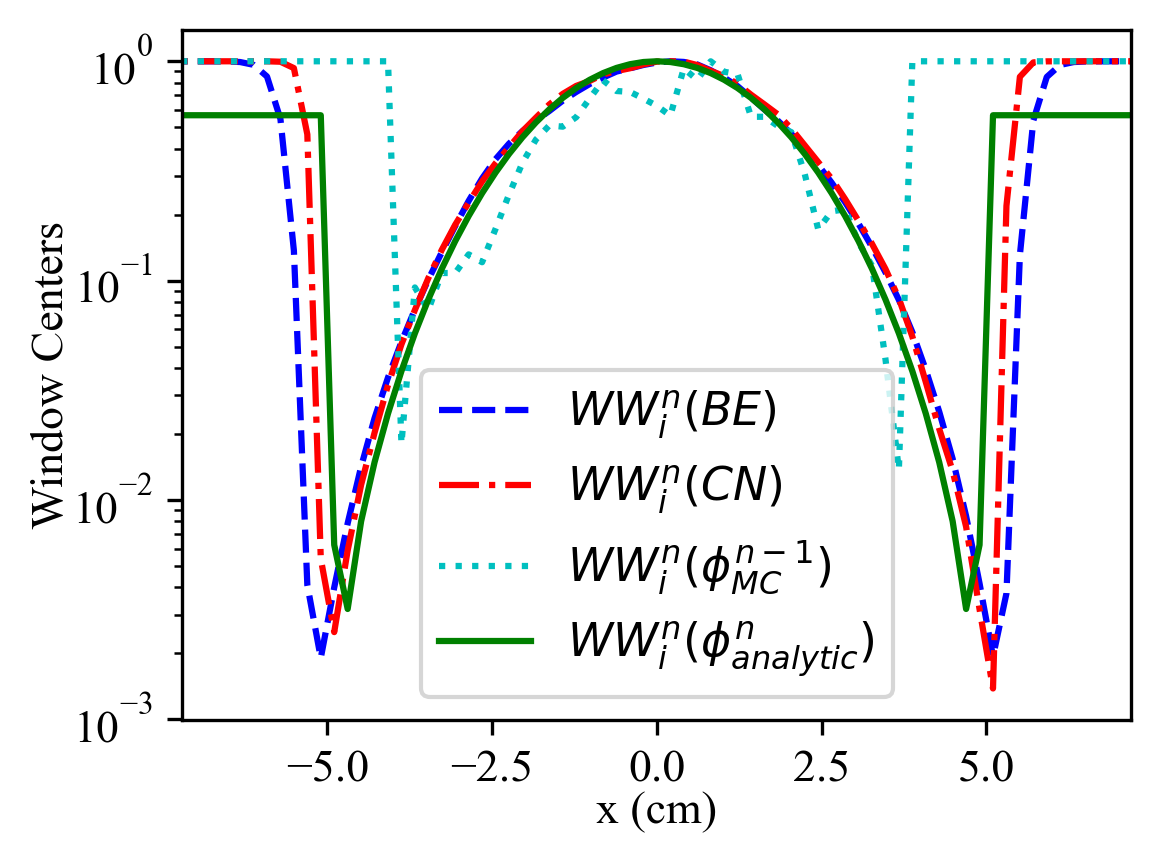}
        \caption{Modified weight window centers for $t~\in~[4.5,5]$}
        \label{fig:mod_centers}
    \end{minipage}\hfill
    \begin{minipage}{0.5\textwidth}
        \centering
        \includegraphics{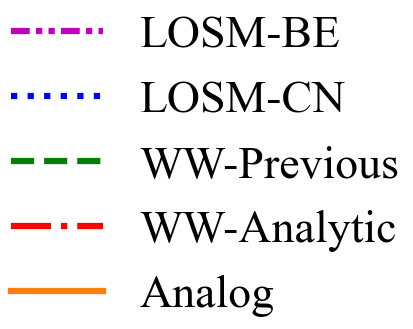}
        \caption{Legend for all future plots (Figs. \ref{fig:azur_flux_3}-\ref{fig:azur_tilde_flux_-1})}
        \label{fig:legend}
    \end{minipage}\hfill
    \begin{minipage}{0.5\textwidth}
        \centering
        \includegraphics[width=\linewidth]{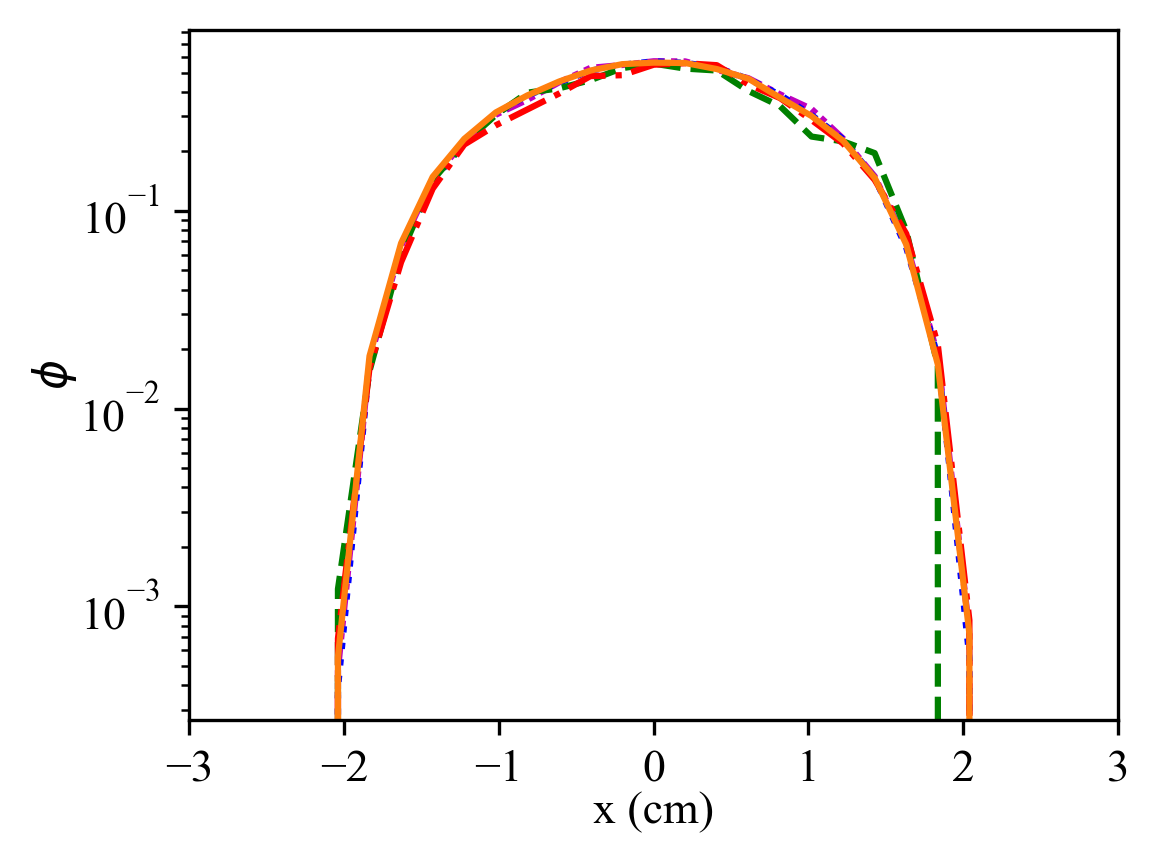}
        \caption{Scalar flux for $t\in[1.5,2]$ s}
        \label{fig:azur_flux_3}
    \end{minipage}\hfill
    \begin{minipage}{0.5\textwidth}
        \centering
        \includegraphics[width=\linewidth]{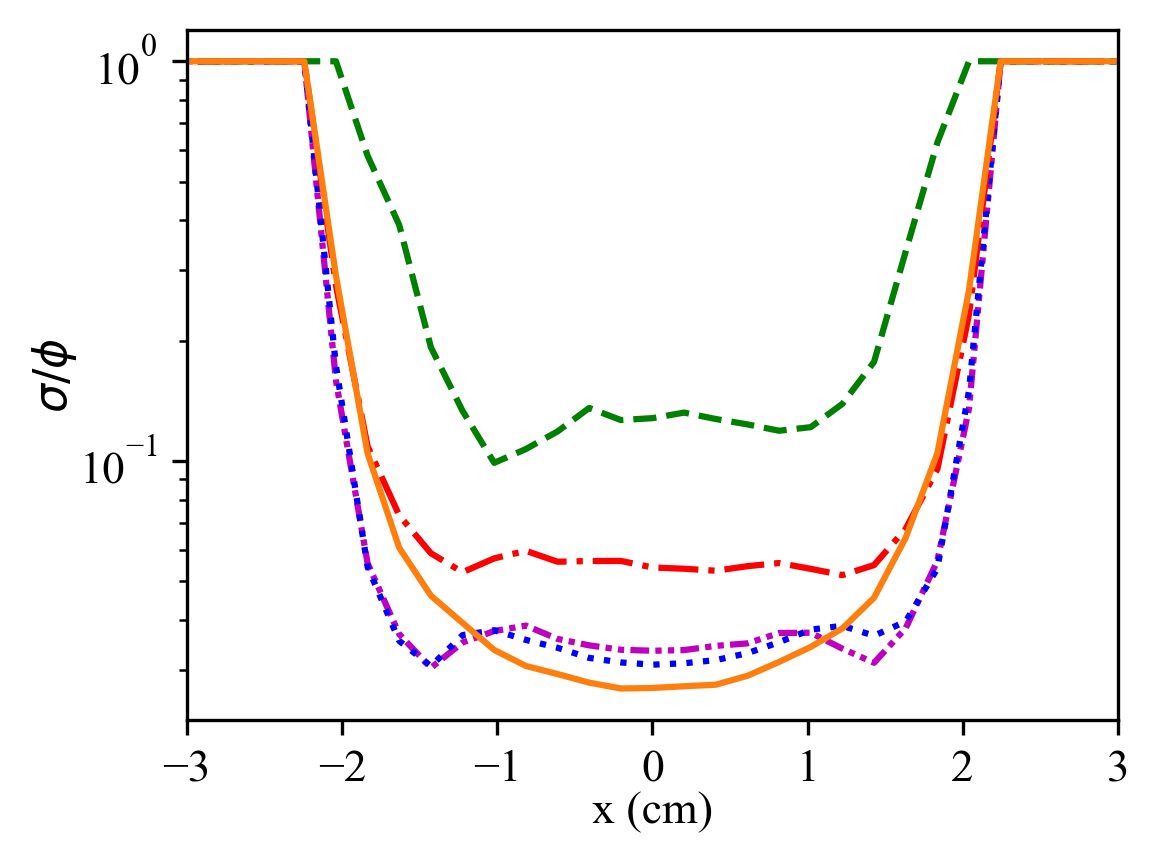}
        \caption{Relative standard deviation for $t~\in~[1.5,2]$ s}
        \label{fig:azur_rel_sdev_3}
    \end{minipage}\hfill
    
        \begin{minipage}{0.5\textwidth}
    	\centering
    	\includegraphics[width=\linewidth]{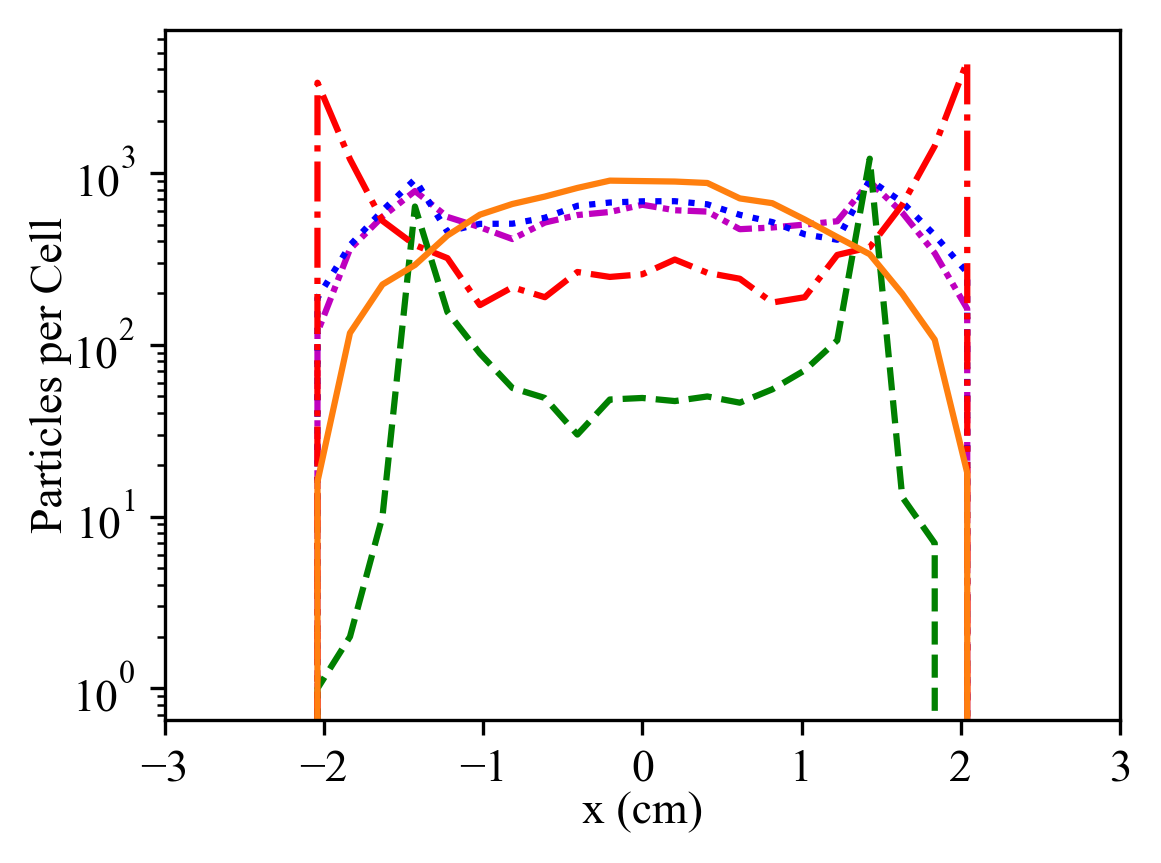}
    	\caption{Particle distribution at $t = 1.5$ s}
    	\label{fig:azur_particle_dist_3}
    \end{minipage}\hfill
    \begin{minipage}{0.5\textwidth}
    	\centering
    	\includegraphics[width=\linewidth]{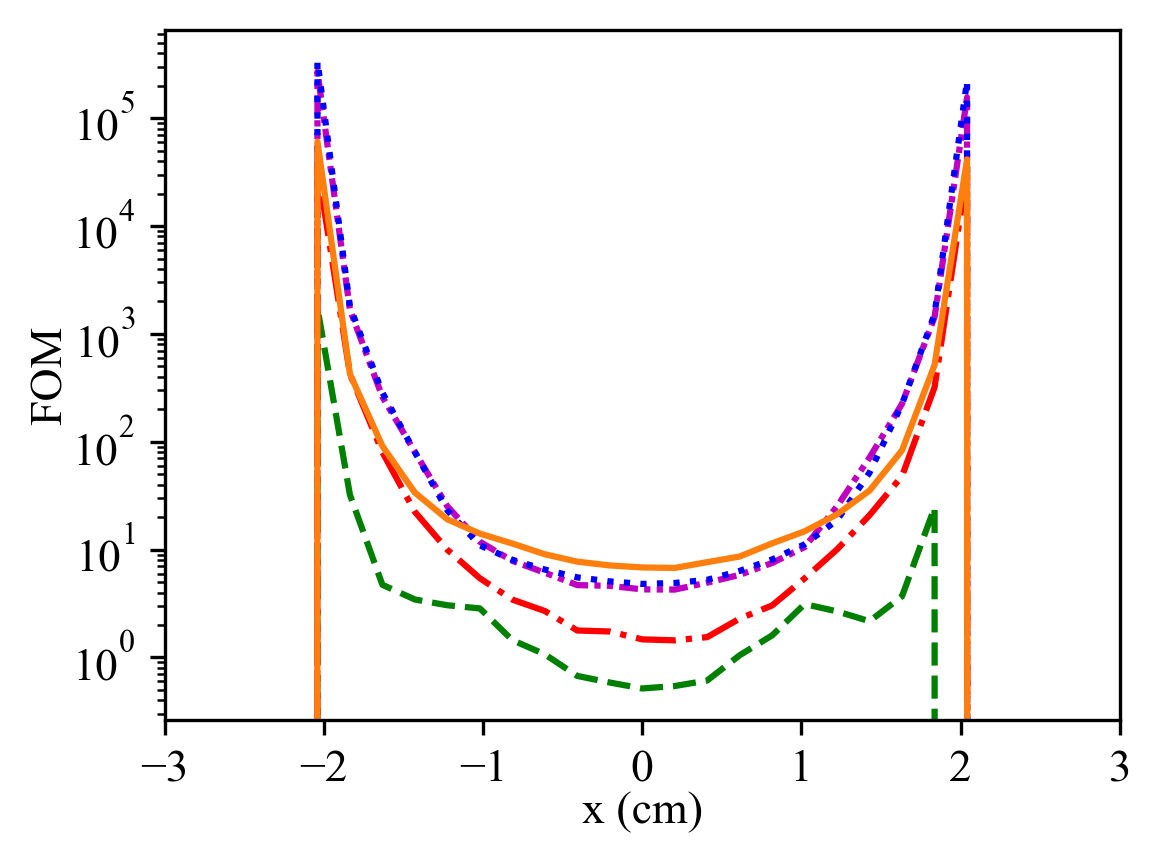}
    	\caption{Spatial FOM for $t\in[1.5,2]$ s}
    	\label{fig:azur_FOM_3}
    \end{minipage}\hfill
    
\end{figure}

\begin{figure}[t!]
    \centering

    	\begin{minipage}{0.5\textwidth}
	\centering
	\includegraphics[width=\linewidth]{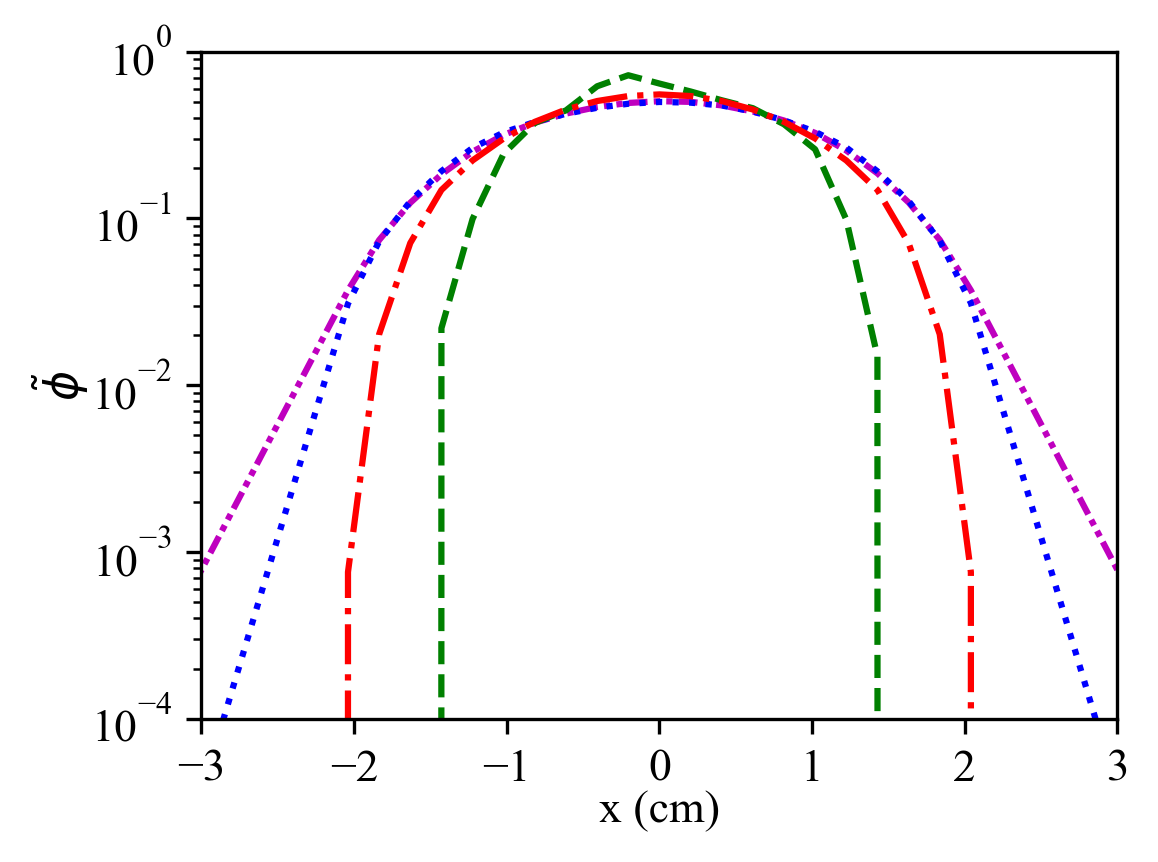}
	\caption{Auxiliary solution   $\tilde \phi^n$ at $t = 2$ s}
	\label{fig:azur_tilde_flux_3}
\end{minipage}\hfill
\begin{minipage}{0.5\textwidth}
	\centering
	\includegraphics[width=\linewidth]{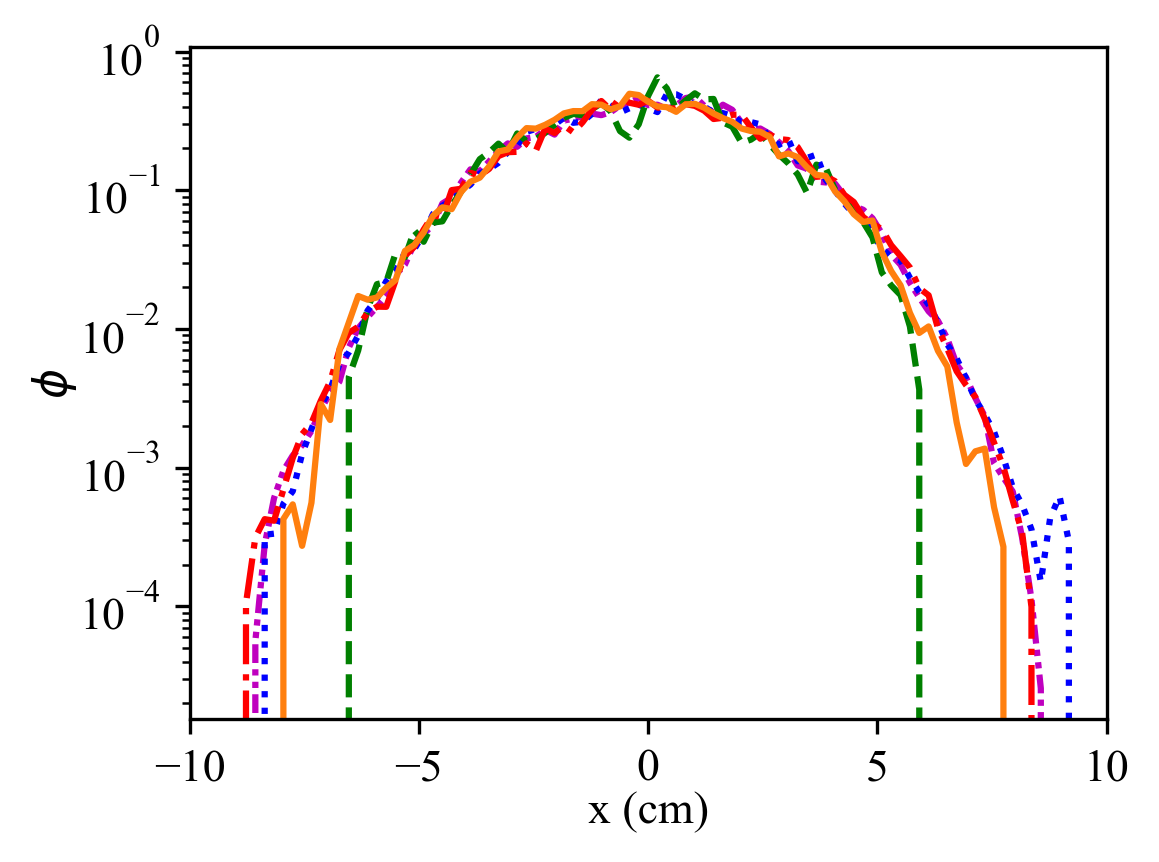}
	\caption{Scalar flux for $t\in[9.5,10]$ s}
	\label{fig:azur_flux_-1}
\end{minipage}\hfill

    \begin{minipage}{0.5\textwidth}
        \centering
        \includegraphics[width=\linewidth]{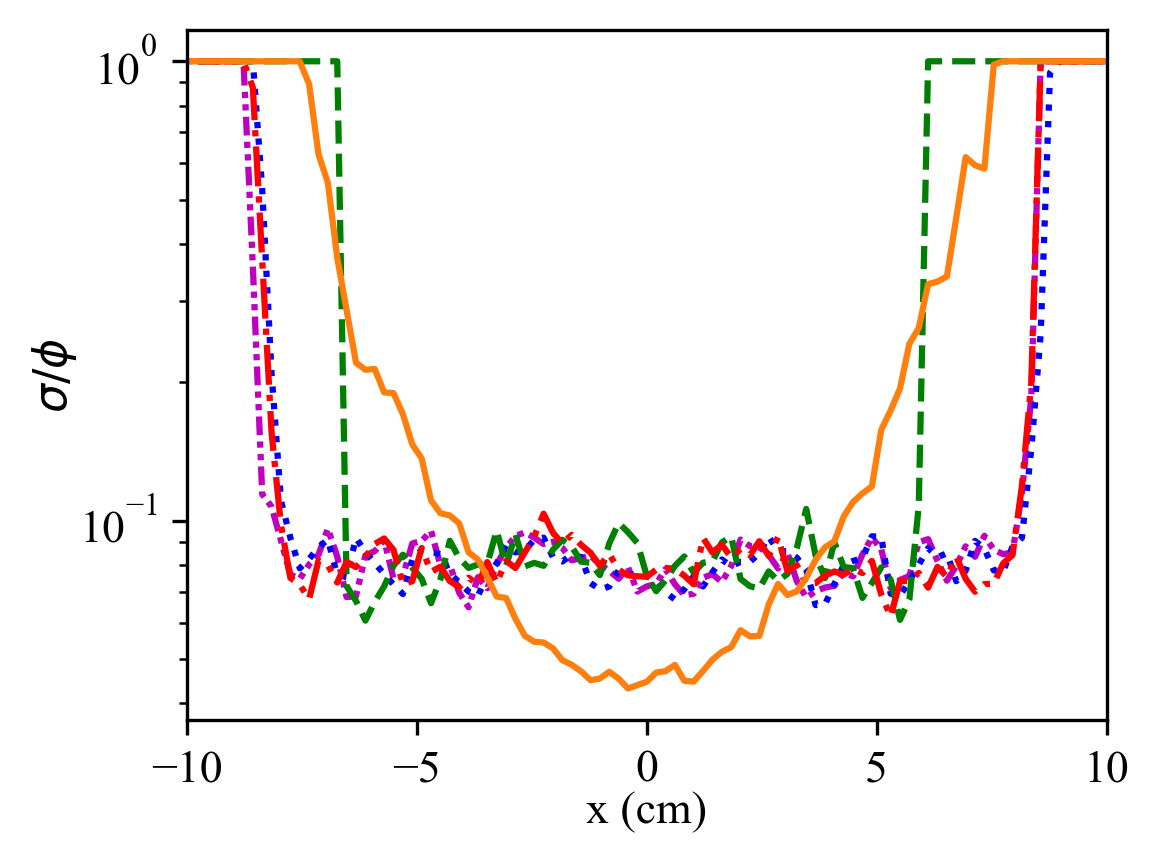}
        \caption{Relative standard deviation for $t~\in~[9.5,10]$ s}
        \label{fig:azur_rel_sdev_-1}
    \end{minipage}\hfill
    \begin{minipage}{0.5\textwidth}
        \centering
        \includegraphics[width=\linewidth]{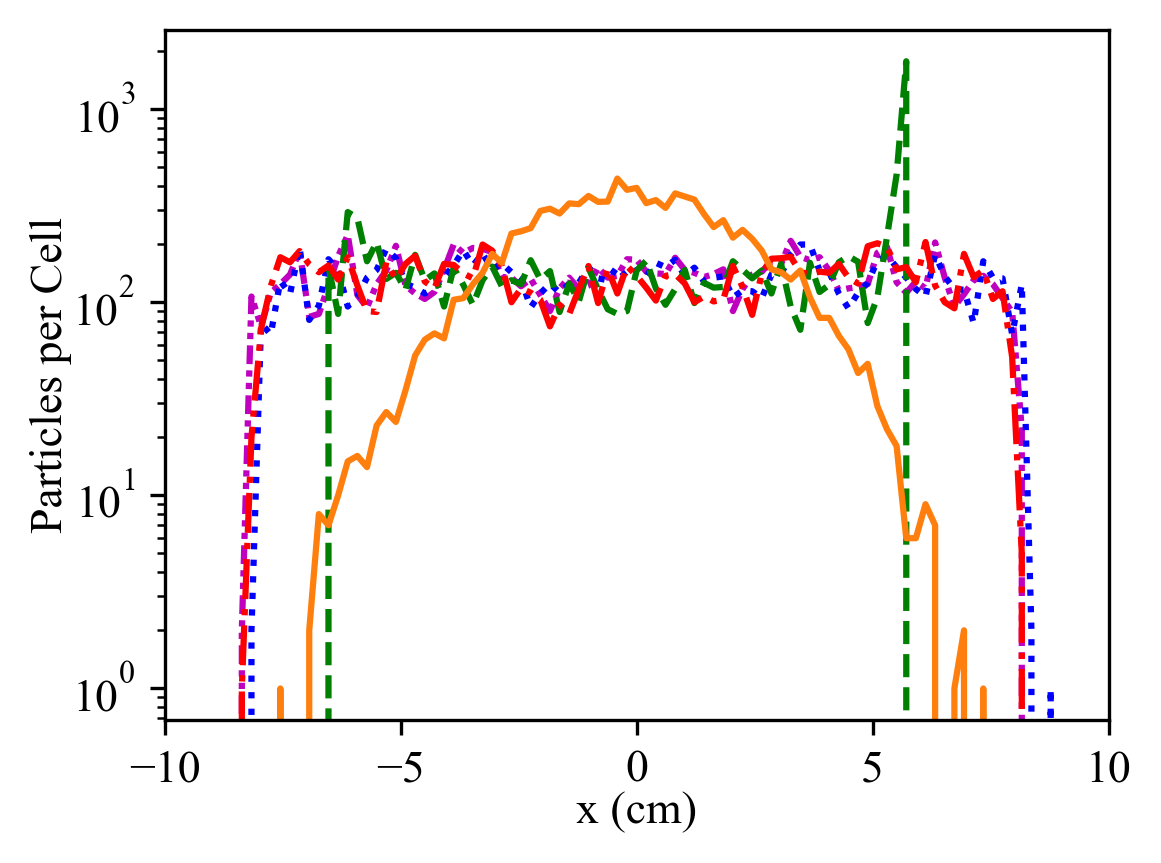}
        \caption{Particle distribution at $t = 9.5$ s\\ \text{}}
        \label{fig:azur_particle_dist_-1}
    \end{minipage}\hfill    
    \begin{minipage}{0.5\textwidth}
        \centering
        \includegraphics[width=\linewidth]{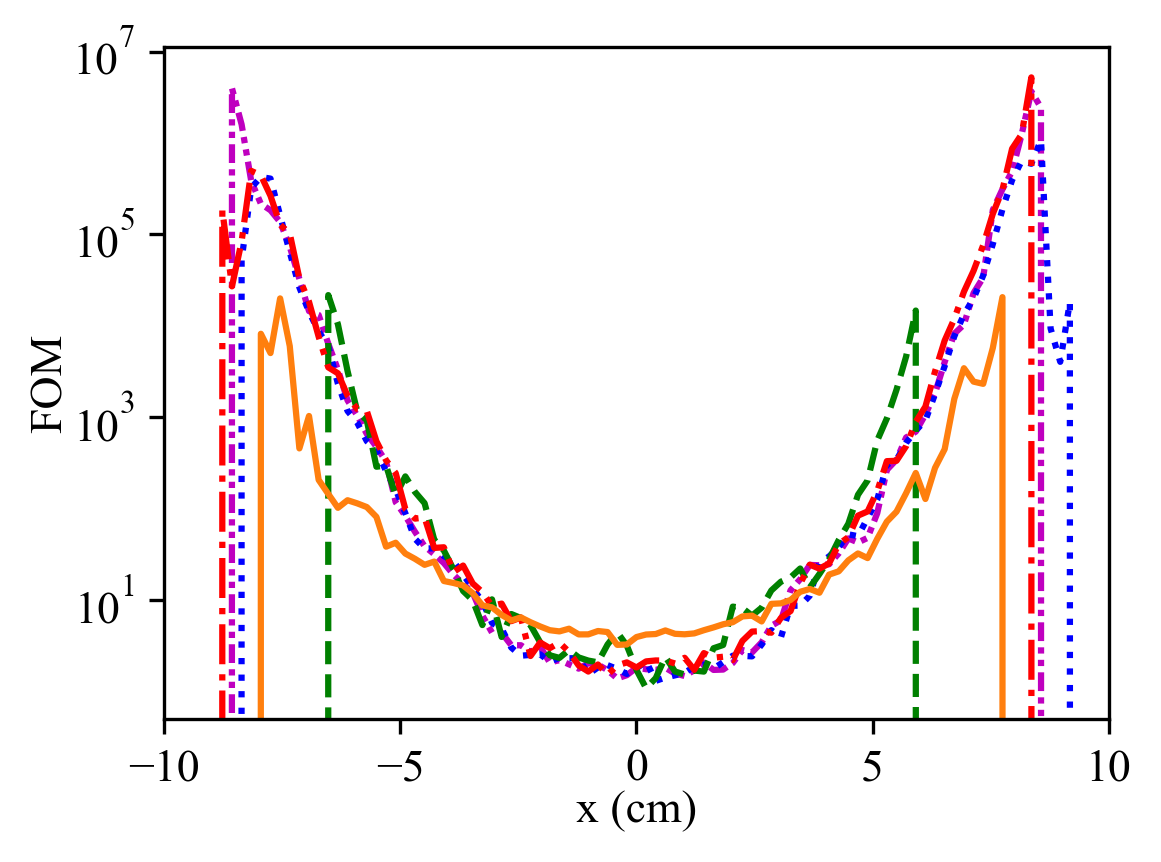}
        \caption{Spatial FOM for $t\in[9.5,10]$ s}
        \label{fig:azur_FOM_-1}
    \end{minipage}\hfill
\centering
    \begin{minipage}{0.5\textwidth}
    \centering
    \includegraphics[width=\linewidth]{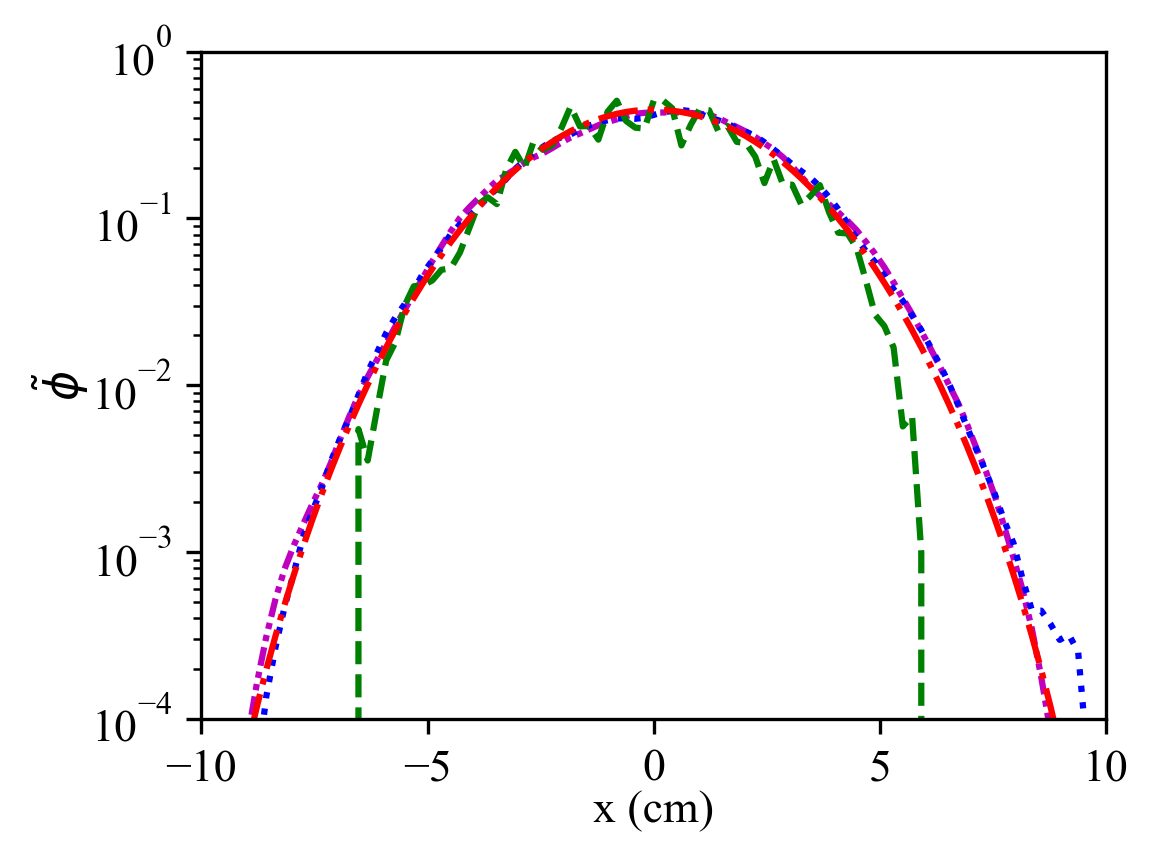}
    \caption{Auxiliary solution  $\tilde \phi^n$ at $t = 10$ s}
    \label{fig:azur_tilde_flux_-1}
\end{minipage}\hfill
\end{figure}

In Figure \ref{fig:mod_centers} the weight window centers are plotted after the modification for waves was applied. The scalar flux in Figures \ref{fig:azur_flux_3} and \ref{fig:azur_flux_-1} was computed using track-length tallies. These results show that WW-Analytic captures the wavefront position most accurately, followed closely by LOSM-CN. The windows computed in WW-Previous restrict particles from reaching the wavefront. This is because the wavefront position in the lagged solution has not traveled as far, and the modification applied to the centers causes rouletting of particles before they reach the true position of the wavefront. The analog solution also does not accurately capture the wavefront due to the low particle counts in that region. 
Note that Figure \ref{fig:legend} shows legend for  plots in Figures \ref{fig:azur_flux_3}-\ref{fig:azur_tilde_flux_-1}.

Figures \ref{fig:azur_FOM_3} and \ref{fig:azur_FOM_-1} show the spatial distribution of the figure of merit (FOM) defined as $(\sigma^2T)^{-1}$, where $T$ is the runtime for the timestep. The LOSM-CN shows improvement in FOM on the outer regions of the problem, at the expense of FOM in the center. This is expected since the effect of the windows can be interpreted as distributing particles from the center to the edges in order to achieve a uniform distribution. This uniformity can be seen in Figures \ref{fig:azur_particle_dist_3} and \ref{fig:azur_particle_dist_-1}, which show the distribution of Monte Carlo particles in space. Having a uniform particle distribution also ensures a uniform relative standard deviation,
\begin{equation}
	\sigma_{relative,i}^n=\frac{\sigma_i^n}{\phi_i^n} 
	\vspace{-0.15cm}
\end{equation}
as seen in Figures \ref{fig:azur_rel_sdev_3} and \ref{fig:azur_rel_sdev_-1}.

In Figures \ref{fig:azur_tilde_flux_3} and \ref{fig:azur_tilde_flux_-1} the auxiliary solutions used for the weight windows are plotted. Here the differences between the WW-Previous, LOSM-CN, LOSM-BE, and WW-Analytic are highlighted. This figure also shows that the CN scheme captures the slope of the wavefront more accurately than the BE scheme, as expected due to the reduced numerical diffusion of the second order scheme. Figures \ref{fig:azur_particle_dist_3} and \ref{fig:azur_particle_dist_-1} indicate that the LOSM-CN weight windows achieve a particle distribution that more closely aligns with the WW-Analytic than WW-Previous does.

\section{Conclusions}\label{sec:conclusion}

In this work, we
have developed a new algorithm for automated  weight windows for Monte Carlo time-dependent particle transport using the solution of the hybrid LOSM  problem for defining weight window centers. The time-dependent LOSM equations are approximated by means of the first- and second-order time-integration schemes. Our results indicate that this method modifies the particle distribution in a way that improves computational efficiency, especially in traveling wave problems at the wavefront. We demonstrated the ability to obtain a more accurate forward hybrid solution using a scheme with updates to parameters of the LOSM problem during a timestep. Future work will expand these hybrid methods to higher dimensions, as well as investigating advanced noise reduction techniques.

\section*{Acknowledgments}

This work was supported by the Center for Exascale Monte-Carlo Neutron Transport (CEMeNT) a PSAAP-III project funded by the Department of Energy,x
 grant number: DE-NA003967.

\bibliographystyle{elsarticle-num}
\bibliography{references}

\end{document}